\documentclass[11 pt]{article}

\usepackage{amssymb}

\def\nmark{\mbox{$\rm\bf\kern0.2em\rule{0.06em}{1.45ex}\kern-0.3em
N$}}
\def\dmark{\mbox{$\rm\bf\kern0.2em\rule{0.06em}{1.45ex}\kern-0.3em
D$}}
\def\cmark{\mbox{$\rm\bf\kern0.2em\rule{0.06em}{1.45ex}\kern-0.3em
C$}}
\def\rmark{\mbox{$\rm\bf\kern0.2em\rule{0.06em}{1.45ex}\kern-0.3em
R$}}

\begin{document}
\title{\large \bf
Weighted composition operators on the Fock space}
 \author{Mahsa Fatehi }

{\maketitle}
\begin{abstract}
In this paper, we study weighted composition operators on the Fock space. We show that a weighted composition operator is cohyponorma if and only if it is normal. Moreover,  we give a complete characterization of  closed range weighted composition operators. Finally, we find norms of some weighted composition operators.
\end{abstract}


\section{Introduction}
Let $\mathbb{D}$ denote the open unit disk in the complex plane $\mathbb{C}$.
 For an analytic map $\varphi$, let $ \varphi_{0}$ be the identity function, $\varphi_{1}=\varphi$ and $\varphi_{n+1} = \varphi \circ \varphi_{n}$ for $n=1,2,...$. We
call them the iterates of $\varphi$. It is well-known that if $\varphi$, neither the identity  nor an elliptic automorphism of $\mathbb{D}$ (i.e., $\varphi$ is an automorphism of  $\mathbb{D}$ with a fixed point in $\mathbb{D}$), is an analytic map on the unit disk into itself, then there exists  a point $w$ in $\overline{\mathbb{D}}$  such that  $\varphi_n$  converges to $w$ uniformly on compact subsets of $\overline{\mathbb{D}}$. The point $w$ is called the  Denjoy-Wolff point of $\varphi$. The  Denjoy-Wolff point $w$ is the unique fixed point of $\varphi$ in $\overline{\mathbb{D}}$ so that $|\varphi'(w)|\leq 1$ (see \cite{cm1}).\par
 Recall that the Fock space $\mathcal{F}^{2}$ is a Hilbert space of all entire  functions on  $\mathbb{C}$ that are square integrable with respect to the Gaussian measures $d \mu(z)=\pi^{-1}e^{-|z|^{2}}dA(z)$, where $dA$ is the usual Lebesgue measure on $\mathbb{C}$. The Fock space $\mathcal{F}^{2}$ is a reproducing kernel Hilbert space with inner product
$$\langle f, g \rangle=\int_{\mathbb{C}}f(z)\overline{g(z)}d\mu(z)$$
and reproducing kernel function $K_{w}(z)=e^{\overline{w}z}$ for any $w \in \mathbb{C}$. Note that for any $w \in \mathbb{C} $, $\|K_{w}\|=e^{|w|^{2}/2}$. For each $w \in \mathbb{C}$, we define the normalized reproducing kernel as $k_{w}(z)=\frac{K_{w}(z)}{\|K_{w}\|}=e^{\overline{w}z-|w|^{2}/2}$. For each nonnegative integer $n$, let $e_{n}(z)=z^{n}/\sqrt{n!}$. The set $\{e_{n}\}$ is an orthonormal basis for $\mathcal{F}^{2}$ (see \cite{z}).\par
For entire function $\varphi$ on $\mathbb{C}$, the composition operator $C_{\varphi}$ on $\mathcal{F}^{2}$ is defined as $C_{\varphi}(f)=f\circ\varphi$ for any $f \in \mathcal{F}^{2}$; moreover, for $\psi \in \mathcal{F}^{2}$, the weighted composition operator $C_{\psi,\varphi}$ is defined by $C_{\psi,\varphi}f=\psi\cdot (f\circ\varphi)$. There exists some literature on composition operators acting on the Hardy and Bergman spaces. The books \cite{cm1} and \cite{sh} are the  important references.

Carswell et al. \cite{car} characterized the bounded and compact  composition operators on the Fock space over $\mathbb{C}^{n}$. Specified to the one-dimensional case, they stated that $C_{\varphi}$ is bounded if and only if $\varphi(z)=az$, where $|a|=1$ or $\varphi(z)=az+b$ with $|a|< 1$.
  In \cite{ueki}, Ueki found the  criteria to characterize the boundedness and compactness of weighted composition operators on the Fock space. After that in \cite{lef}, Le obtained much easier criteria for the boundedness and compactness of weighted composition operators on the Fock space.
On the Hardy space,  normal weighted composition operators were studied; moreover, unitary weighted composition operators were characterized (see \cite{Bourdon}).  Also in \cite{coko}, cohyponormal weighted composition operators were obtained. After that in \cite{hypomahsa}, hyponormal weighted composition operators were investigated on the Hardy and weighted Bergman spaces. Unitary weighted composition operators $C_{\psi,\varphi}$  on the Fock space were obtained in \cite{zhao1}. Also invertible weighted composition operators on the Fock space were characterized in \cite{zhao2}.
In the second section,  we find normaloid, hyponormal and cohyponormal composition operators. After that we obtain all hyponormal weighted composition operators $C_{\psi,\varphi}$, where $\psi=K_{c}$ for each $c \in \mathbb{C}$. Moreover, we study a class of normaloid weighted composition operators. Next, we show that for $\varphi(z)=az+b$, $C_{\psi,\varphi}$ is cohyponormal if and only if $\psi=\psi(0)K_{b\frac{\overline{a}-1}{a-1}}$.
Closed range composition operators were studied on the Hardy and weighted Bergman spaces in  \cite{af}, \cite{gh} and \cite{zo}. In the third section, we characterize closed range weighted composition operators on the Fock space. In the fourth section, we find norm of $C_{\psi,\varphi}$ on $\mathcal{F}^{2}$, when $\psi=K_{c}$ for any $c \in \mathbb{C} $.\\

\section{Normaloid weighted composition operators}
Suppose that $T$ is a bounded operator on a Hilbert space. Throughout this paper, the spectrum of $T$, the essential spectrum of $T$, and the point spectrum of $T$ are denoted by $\sigma(T)$, $\sigma_{e}(T)$,  $\sigma_{p}(T)$ respectively. Also the spectral radius of $T$ is denoted by $r(T)$.\\
Le \cite{lef} studied the boundedness of weighted composition operator on the Fock space. His  result shows that if $C_{\psi,\varphi}$ is bounded on $\mathcal{F}^{2}$, then $\varphi(z)=az+b$, where $|a|\leq 1$; furthermore, he proved that if $|a|=1$ and $C_{\psi,az+b}$ is bounded on $\mathcal{F}^{2}$, then $\psi=\psi(0)K_{-\overline{a}b}$. We use these facts frequently in this paper and so  throughout this paper, we assume that  $\varphi(z)=az+b$, where $|a|\leq 1$.
\cite{car}, \cite{zhao1} and \cite{zhao} were written on another Fock space (see \cite[p. 33]{z}), but their results hold for $\mathcal{F}^{2}$ which is considered in this paper, with identical arguments. \par
In the following proposition, we investigate $\sigma_{p}(C_{\psi,\varphi})$, when $\varphi(z)=az+b$ and $|a|< 1$. Note that in the case that $|a|=1$, as we saw in the preceding paragraph, $\psi=\psi(0)K_{-\overline{a}b}$. Then  $C_{\psi,\varphi}$ is a constant multiple of a unitary operator from  \cite[Corollary 1.2]{zhao1}. Moreover, in \cite[Corollary 1.4]{zhao1}, the spectrum of unitary weighted composition operators were characterized.
\\ \par

{\bf Proposition 2.1.} {\it Let $\psi$ and $\varphi$ be entire functions on $\mathbb{C}$ and $C_{\psi,\varphi}$ be a bounded weighted composition operator on $\mathcal{F}^{2}$. Suppose that $\varphi(z)=az+b$, where $|a|< 1$ and $b \in \mathbb{C}$. If $\lambda \in \sigma_{p}(C_{\psi,\varphi})$, then $|\lambda|\leq \psi(\frac{b}{1-a})$. Moreover, if $\psi(\frac{b}{1-a})=0$ and $\varphi$ and $\psi$ are nonconstant, then $C_{\psi,\varphi}$ has no eigenvalues.}\bigskip

{\bf Proof.} Suppose that $|a|< 1$. First we find a weighted composition operator $C_{\widetilde{\psi},\widetilde{\varphi}}$ which is unitary equivalent to $C_{\psi,\varphi}$ such that the fixed point of $\widetilde{\varphi}$ lies in $\mathbb{D}$ and $\widetilde{\varphi}$ is a self-map of $\mathbb{D}$. There exists a positive integer $N$ such that $|a|+\frac{1}{N}|1-a|< 1$. It is not hard to see that there is a complex number
$u \in\mathbb{C}$ such that $|u+\frac{b}{1-a}|< \frac{1}{N}$. By \cite[Corollary 1.2]{zhao1}, we know that $C_{k_{u},z-u}$ is unitary (the operator  $C_{k_{u},z-u}$ is known as the Weyl unitary). By  \cite[Proposition 3.1]{lef}, $C_{k_{u},z-u}^{\ast}=C_{k_{-u},z+u}$, so

\begin{eqnarray}
C_{k_{u},z-u}C_{\psi,\varphi}C_{k_{-u},z+u}&=&
\frac{1}{\|K_{u}\|^{2}}C_{e^{\overline{u}z},z-u}C_{\psi,\varphi}C_{e^{-\overline{u}z},z+u}\nonumber\\
&=&\frac{1}{\|K_{u}\|^{2}}e^{\overline{u}z}\cdot \psi(z-u)\cdot e^{(-\overline{u}(az+b))\circ(z-u)}C_{(z+u)\circ (az+b)\circ (z-u)}\nonumber\\
&=&\frac{1}{\|K_{u}\|^{2}}e^{\overline{u}z}\cdot \psi(z-u)\cdot e^{-\overline{u}(az-au+b)}C_{az+u(1-a)+b}.
\end{eqnarray}
Let $\widetilde{\varphi}(z)=az+u(1-a)+b$ and
\begin{eqnarray}
\widetilde{\psi}(z)&=&\frac{1}{\|K_{u}\|^{2}}e^{\overline{u}z}\cdot \psi(z-u)\cdot e^{-\overline{u}(az-au+b)}.
\end{eqnarray}
It is easy to see that the fixed point of $\widetilde{\varphi}$ is $u+\frac{b}{1-a}$ which belongs to $\mathbb{D}$. Because $|a|+|u(1-a)+b|< 1$, $\widetilde{\varphi}$ is a self-map of $\mathbb{D}$. Since $C_{\psi,\varphi}$ is unitary equivalent to $C_{\widetilde{\psi},\widetilde{\varphi}}$, $\sigma_{p}(C_{\psi,\varphi})=\sigma_{p}(C_{\widetilde{\psi},\widetilde{\varphi}})$. Assume that $\lambda$ is a nonzero eigenvalue for $C_{\widetilde{\psi},\widetilde{\varphi}}$ with corresponding eigenvector $h$. We obtain
\begin{equation}
\lambda^{n}h(z)=\Pi_{j=0}^{n-1}\widetilde{\psi}(\widetilde{\varphi}_{j}(z))h(\widetilde{\varphi}_{n}(z))
\end{equation}
for each $z\in \mathbb{C}$ and positive integer $n$. For any fixed point $z\in \mathbb{C}$, we obtain
\begin{eqnarray}
|h(\widetilde{\varphi}_{n}(z))|&=&|\langle h\circ\widetilde{\varphi}_{n},K_{z}\rangle|\nonumber\\
&\leq&\|h\circ \widetilde{\varphi}_{n}\|\|K_{z}\|\nonumber\\
&=&\|h\circ\widetilde{\varphi}_{n}\|e^{\frac{|z|^{2}}{2}}\nonumber\\
&=&\|C_{\widetilde{\varphi}_{n}}(h)\|e^{\frac{|z|^{2}}{2}}\nonumber\\
&\leq&\|C_{\widetilde{\varphi}_{n}}\|\|h\|e^{\frac{|z|^{2}}{2}}.
\end{eqnarray}
Since $h$ is not the zero function, we can choose $z \in \mathbb{D}$ such that $h(z) \neq 0$. Since $u+\frac{b}{1-a}$ is the Denjoy-Wolff point of $\widetilde{\varphi}$, $\widetilde{\varphi}_{j}(z)\rightarrow u+\frac{b}{1-a}$ and $\widetilde{\psi}(\widetilde{\varphi}_{j}(z))\rightarrow \widetilde{\psi}(u+\frac{b}{1-a})$ as $j\rightarrow \infty$. Take $n$-th roots of the absolute value each side of Equation (3), use Equation (4) and let $n\rightarrow \infty$, we get $|\lambda|\leq |\widetilde{\psi}(u+\frac{b}{1-a})|r(C_{\varphi})$. From Equation (2) and \cite[Theorem 1.1]{zhao}, we see that $|\lambda|\leq |\widetilde{\psi}(u+\frac{b}{1-a})|r(C_{\varphi})=|\psi(\frac{b}{1-a})|$. Then
\begin{equation}
|\lambda|\leq |\psi(\frac{b}{1-a})|.
\end{equation}
Now assume that $\widetilde{\psi}(u+\frac{b}{1-a})=\psi(\frac{b}{1-a})=0$ and $\varphi$ and $\psi$ are nonconstant. Thus, by Equation (5), $\lambda=0$ is the only possible eigenvalue for $C_{\psi,\varphi}$. Since $\psi$ is not the zero function and $\varphi$ is not  constant, the Open Mapping Theorem implies that $0$ cannot be an eigenvalue for $C_{\psi,\varphi}$ (see  idea of the proof of \cite[Lemma 4.1]{Bourdon1}).\hfill $\Box$ \\ \par

Suppose that $T$ is a bounded operator. The  operator $T$ is hyponormal (cohyponormal) if $T^{\ast}T\geq TT^{\ast}$ ($T^{\ast}T\leq TT^{\ast}$). Also $T$ is normaloid if $\|T\|=r(T)$. It is well known that hyponormal (cohyponormal)  operators are normaloid. In the following proposition, we characterize hyponormal, cohyponormal and normaloid composition operators on $\mathcal{F}^{2}$.\\ \par

{\bf Proposition 2.2.} {\it Let $\varphi(z)=az+b$, where $|a|\leq1$ and $b \in \mathbb{C}$. Then $C_{\varphi}$ is a bounded normaloid (hyponormal or cohyponormal) operator if and only if $b=0$.}\bigskip

{\bf Proof.} Let $C_{\varphi}$ be normaloid (hyponormal or cohyponormal). Suppose that $b \neq 0$. Then \cite[Theorem 1]{car} states that $|a|< 1$. By \cite[Theorem 4]{car}, $\|C_{\varphi}\|=e^{\frac{1}{2}\frac{|b|^{2}}{1-|a|^{2}}}$ (note that the inner product for $\mathcal{F}^{2}$ in this paper is different from  \cite{car}, so the norm of a composition operator is not exactly the same as \cite[Theorem 4]{car}; furthermore, in Remark 4.1, $\|C_{\varphi}\|$ will be described). Also \cite[Theorem 1.1]{zhao} implies that $r(C_{\varphi})=1$. Since $C_{\varphi}$ is normaloid (hyponormal or cohyponormal), $e^{\frac{1}{2}\frac{|b|^{2}}{1-|a|^{2}}}=1$. Hence $b=0$ which is a contradiction. \\
Conversely,  suppose that $\varphi(z)=az$, where $|a|\leq 1$. Invoking \cite[Lemma 2]{car}, $C_{\varphi}^{\ast}=C_{\overline{a}z}$. Then $C_{az}$  is normal and the result follows.\hfill $\Box$ \\ \par

Suppose that $C_{\psi,\varphi}$ is a bounded weighted composition operator and $\varphi(z)=az+b$. Note that if $a=1$, then from \cite[Corollary 1.2]{zhao1} and as we saw in the second paragraph of this section, $C_{\psi,\varphi}$ is a constant multiplie of a unitary operator. Thus, $C_{\psi,\varphi}$ is normal, normaloid, hyponormal and cohyponormal. Hence we state the following proposition for $a\neq 1$.\\ \par

{\bf Proposition 2.3.} {\it Let $\psi=K_{c}$ and $\varphi(z)=az+b$, where $|a| \leq 1$, $a\neq 1$ and $b,c \in \mathbb{C}$. Suppose that $C_{\psi,\varphi}$  is bounded on $\mathcal{F}^{2}$. Then the following are equivalent.\\
(a) $C_{\psi,\varphi}$ is hyponormal.\\
(b) $C_{\psi,\varphi}$ is cohyponormal.\\
(c) $C_{\psi,\varphi}$ is normaloid.\\
(d) $c=b \frac{\overline{a}-1}{a-1}$.}\bigskip

{\bf Proof.} There is $u \in \mathbb{C}$ such that $c=u(\overline{a}-1)$. By Equation (1) and some calculation, $C_{\psi,\varphi}$ is unitarily equivalent to $C_{\widetilde{\psi},\widetilde{\varphi}}$, where $\widetilde{\psi}(z)=e^{-|u|^{2}}e^{\overline{u}z}\cdot \psi(z-u)\cdot e^{-\overline{u}(az-au+b)}=\psi(\frac{b}{1-a})$
and $\widetilde{\varphi}(z)=az+u(1-a)+b$. We can see that $C_{\psi,\varphi}$ is unitarily equivalent to $\psi(\frac{b}{1-a})C_{az+u(1-a)+b}$.\par
(a) $\Rightarrow$ (d). Suppose that $C_{\psi,\varphi}$ is hyponormal. Then $C_{az+u(1-a)+b}$ is hyponormal. Proposition 2.2 implies that $u(1-a)+b=0$. Since $u=\frac{c}{\overline{a}-1}$, we conclude that $c=b \frac{\overline{a}-1}{a-1}$.\par
(d) $\Rightarrow$ (a). Assume that $c=b \frac{\overline{a}-1}{a-1}$. Let $u=\frac{b}{a-1}$. By Equation (1) and some calculation, $C_{\psi,\varphi}$ is unitarily equivalent to $e^{\frac{|b|^{2}}{1-\overline{a}}}C_{\widetilde{\varphi}}$, where $\widetilde{\varphi}(z)=az$. We infer from Proposition 2.2  that  $C_{\widetilde{\varphi}}$ is hyponormal and so $C_{\psi,\varphi}$
is hyponormal. \par
(b) $\Leftrightarrow$ (d). The idea of the proof  is similar to (a) $\Leftrightarrow$ (d).\par
(c) $\Leftrightarrow$ (d). The idea of the proof  is similar to (a) $\Leftrightarrow$ (d).\hfill $\Box$ \\ \par
%(b) $\rightarrow$ (d). Suppose that $C_{\psi,\varphi}$ is cohyponormal. Then $C_{\psi,\varphi}$ is normailoid and so the statement (d)
%holds.\par
%(d) $\rightarrow$ (b). Suppose that the statement (d) occurs. Then $C_{\psi,\varphi}$ is normal by \cite[Theorem 3.3]{lef}. It shows that $C_{\psi,\varphi}$ is %cohyponormal.\hfill $\Box$ \\ \par

Suppose that $\varphi(z)=az+b$, where $|a|=1$ and $\psi$ is an entire function. If $C_{\psi,\varphi}$ is bounded on $\mathcal{F}^{2}$, then as we saw in the second paragraph of this section, $\psi(z)=\psi(0)K_{-\beta}$, where $\beta=\overline{a}b$. Hence $C_{\psi,\varphi}$ is a constant multiple of the unitary operator (see \cite[Corollary 1.2]{zhao1}). It shows that in this case $C_{\psi,\varphi}$ is normaloid and so in the next theorem, we assume that $|a|<1$. \\ \par

{\bf Theorem 2.4.} {\it Suppose that $\psi$ is an entire function and is not identically zero. Assume that $\varphi(z)=az+b$, where $|a|<1$ and $b \in \mathbb{C}$. Let for each $\lambda \in \sigma_{e}(C_{\psi,\varphi})$, $|\lambda| \leq |\psi(\frac{b}{1-a})|$. Then $C_{\psi,\varphi}$ is normaloid if and only if $\psi=\psi(0) K_{b\frac{\overline{a}-1}{a-1}}$.}\bigskip

 {\bf Proof.} Let $p=\frac{b}{1-a}$ (note that it is obvious that $p$ is the  fixed point of $\varphi$). By  \cite[Corollary 1.2]{zhao1}, $C_{k_{p},z-p}$  is unitary. Furthermore, \cite[Proposition 3.1]{lef} states that $C_{k_{p},z-p}^{\ast}=C_{k_{-p},z+p}$.  We obtain
\begin{eqnarray}
H:&=&C_{k_{p},z-p}^{\ast}C_{\psi,\varphi}C_{k_{p},z-p}\nonumber\\
&=&C_{k_{-p},z+p}C_{\psi,\varphi}C_{k_{p},z-p}\nonumber\\
&=&C_{q,\widetilde{\varphi}},
\end{eqnarray}
where
\begin{eqnarray}
\widetilde{\varphi}(z)&=&\varphi(z+p)-p=a(z+p)+b-p=az
\end{eqnarray}
and
\begin{eqnarray}
q(z)&=&k_{-p}(z)k_{p}(\varphi(z+p))\psi(z+p)\nonumber\\
&=&e^{\overline{p}b+|p|^{2}(a-1)}e^{\overline{p}(a-1)z}\psi(z+p)\nonumber\\
&=&e^{\overline{p}(a-1)z}\psi(z+p).
\end{eqnarray}
Since $C_{\psi,\varphi}$ is unitary equivalent to $C_{q,\widetilde{\varphi}}$, $\sigma_{e}(C_{\psi,\varphi})=\sigma_{e}(C_{q,\widetilde{\varphi}})$. It is not hard to see that $q(0)=\psi(p)$. Now suppose that $C_{\psi,\varphi}$ is normaloid. Then $C_{q,\widetilde{\varphi}}$ is normaloid. Since $C_{\psi,\varphi}$ and $C_{q,\widetilde{\varphi}}$ are unitary equivalent, for each $\lambda \in \sigma_{p}(C_{q,\widetilde{\varphi}})$, $|\lambda|\leq |q(0)|$.
We infer from \cite[Proposition 6.7, p. 210]{c1} and \cite[Proposition 4.4, p. 359]{c1} that $r(C_{q,\widetilde{\varphi}})\leq |q(0)|$. Since $C_{q,\widetilde{\varphi}}$ is normaloid, $|q(0)|\geq \|C_{q,\widetilde{\varphi}}\|\geq \|C_{q,\widetilde{\varphi}}(1)\|=\|q\|$. We know that $\{\frac{z^{m}}{\sqrt{m!}}: m\geq 0\}$ is an orthonormal basis for $\mathcal{F}^{2}$. Then $\|q\|\geq |q(0)|$. It shows that $q$ must be constant. Thus, Equation (8) shows that $\psi(z)\cdot e^{\overline{p}(a-1)(z-p)}$ is  constant. Then $\psi(z)=\psi(0)e^{-\overline{p}(a-1)z}=\psi(0)e^{\overline{b}\frac{a-1}{\overline{a}-1}z}=\psi(0)K_{b\frac{\overline{a}-1}{a-1}}(z)$.\\
Conversely, suppose that $\psi=\psi(0)K_{b\frac{\overline{a}-1}{a-1}}$. By Equation (8), $q$ is  constant. Equations (6) and (7) state that $C_{\psi,\varphi}$ is unitarily equivalent to a constant multiple of $C_{az}$. Since by Proposition 2.2, $C_{az}$ is normaloid, $C_{\psi,\varphi}$ is also normaloid.\hfill $\Box$ \\ \par

 We know that for each $c \in \mathbb{C}$, $C_{K_{c},\varphi}$ is compact, where $\varphi(z)=az+b$ and $|a|<1$ (see \cite[Corollary 2.4 ]{zhao}). Hence $\sigma_{e}(C_{K_{c},\varphi})=\{0\}$. Thus,  $C_{K_{c},\varphi}$ satisfies the conditions of Theorem 2.4 and so if $C_{K_{c},\varphi}$ is normaloid, then $c$ must be $b\frac{\overline{a}-1}{a-1}$ (see also Proposition 2.3).\\ \par

{\bf Corollary 2.5.} {\it Suppose that $\psi$ is an entire function and is not identically zero. Assume that $\varphi(z)=az+b$, where $|a|<1$ and $b \in \mathbb{C}$.  Then $C_{\psi,\varphi}$ is  compact and normaloid if and only if $\psi=\psi(0) K_{b\frac{\overline{a}-1}{a-1}}$.}\bigskip

Suppose that $\psi$ is an entire function and $\varphi(z)=az+b$, where $|a|\leq 1$. In the next theorem, we show that $C_{\psi,\varphi}$ is cohyponormal if and only if $C_{\psi,\varphi}$ is normal (see \cite[Theorem 3.3]{lef} and note that if $|a|=1$ and $a\neq 1$, then $\frac{\overline{a}-1}{a-1}=-\overline{a}$).\\ \par

{\bf Theorem 2.6.} {\it Suppose that $\psi$ is an entire function and $\varphi(z)=az+b$, where $|a|\leq 1$. Then $C_{\psi,\varphi}$ is cohyponormal if and only if $\psi=\psi(0) K_{b}\frac{\overline{a}-1}{a-1}$ for $a \neq 1$ and $\psi=\psi(0)K_{-b}$ for $a=1$.}\bigskip

 {\bf Proof.}  We break the proof into two cases. First assume that $a=1$. If $\psi=\psi(0)K_{-b}$, then by \cite[Corollary 1.2]{zhao1}, $C_{\psi,\varphi}$ is a constant multiple of a unitary operator. Thus, $C_{\psi,\varphi}$ is cohyponormal. Now let $C_{\psi,\varphi}$ be cohyponormal. As we saw in the second paragraph of this section $\psi=\psi(0)K_{-b}$.\par
 Now assume that $a\neq 1$. Suppose that $C_{\psi,\varphi}$ is cohyponormal. Let $C_{q,\widetilde{\varphi}}$ be as in Equation (6), where $q$ and $\widetilde{\varphi}$ were obtained in Equations (7) and (8). It is obvious that $C_{q,\widetilde{\varphi}}$ is also cohyponormal. Then $\|C_{q,\widetilde{\varphi}}^{\ast}K_{0}\|  \geq \|C_{q,\widetilde{\varphi}}K_{0}\|$. We obtain $|q(0)|\geq \|q\|$. As we saw in the proof of Theorem 2.4,  $q$ must be  constant, so Equation (8) stats that $\psi(z)=\psi(0)e^{\overline{b}\frac{a-1}{\overline{a}-1}z}$.\\
 The other direction  follows easily from Proposition 2.3. \hfill $\Box$ \\ \par

\section{Closed range weighted composition operator}

In this section, we prove that  $C_{\psi,\varphi}$ has closed range if and only if $C_{\psi,\varphi}$ is a constant multiple of a unitary operator (see \cite[Corollary 1.2 ]{zhao1}).\\ \par

{\bf Theorem 3.1.} {\it Let $\varphi$ and $\psi$ be  entire functions on $\mathbb{C}$ such that $\psi$ is not identically zero.  Suppose that $C_{\psi,\varphi}$ is bounded on $\mathcal{F}^{2}$. Then $C_{\psi,\varphi}$ has closed range if and only if  and $\varphi(z)=az+b$, with $|a| = 1$, $b \in \mathbb{C}$ and $\psi=\psi(0)K_{-\overline{a}b}$.}\bigskip

{\bf Proof.} First suppose that $|a|=1$ and $\psi=(0)K_{-\overline{a}b}$. By \cite[Corollary 1.2 ]{zhao1}, we have $C_{\psi,\varphi}$ is a constant multiple of a unitary operator.  Therefore, $C_{\psi,\varphi}$ has closed range. \\
Conversely, let $C_{\psi,\varphi}$ have closed range  on $\mathcal{F}^{2}$. As we stated in the second paragraph of Section 2,  $\varphi(z)=az+b$, with $|a| \leq 1$.
Suppose that $|a|< 1$. By Equations (6) and (7), $C_{\psi,\varphi}$ is unitarily equivalent to $C_{q,\widetilde{\varphi}}$, where $\widetilde{\varphi}(z)=az$, so without loss of generality, we assume that $\varphi(z)=az$ (note that $C_{\psi,\varphi}$ has closed range if and only if $C_{q,\widetilde{\varphi}}$ has closed range). Since $C_{\psi,\varphi}$ is bounded on $\mathcal{F}^{2}$, $C_{\psi,az}(e^{z})=\psi(z)e^{az}$ belongs to $\mathcal{F}^{2}$. Now we define a bounded linear functional $F_{\psi(z)e^{az}}$ by $F_{\psi(z)e^{az}}(f)=\langle f(z),\psi(z) e^{az} \rangle$ for each $f \in \mathcal{F}^{2}$. We know that $\frac{K_{w}}{\|K_{w}\|}$ converges to zero weakly as $|w|\rightarrow \infty$. Then
$$\lim_{|w|\rightarrow \infty}\langle \frac{K_{w}}{\|K_{w}\|}, \psi(z) e^{az} \rangle=0.$$
It shows that
\begin{eqnarray}
\lim_{|w|\rightarrow \infty}\frac{|\psi(w)||e^{aw}|}{e^{|w|^{2/2}}}&=&0.
\end{eqnarray}
Now if  $a=|a|e^{i\theta}$, we take $w=re^{-i \theta}$, where $r$ is a positive real number. Then $|e^{aw}|=e^{|aw|}$. Equation (9) shows that
\begin{eqnarray}
\lim_{r\rightarrow \infty}|\psi(re^{-i \theta})|^{2}e^{|ar|^{2}-r^{2}}&=0.&
\end{eqnarray}
From Equation (10), we obtain
\begin{eqnarray}
\lim_{r\rightarrow \infty}\|C_{\psi,\varphi}^{\ast}\frac{K_{re^{-i \theta}}}{\|K_{re^{-i \theta}}\|}\|^{2}&=& \lim_{r\rightarrow \infty}|\psi(re^{-i \theta})|^{2}\frac{\|e^{re^{i \theta}\overline{a}z}\|^{2}}{e^{r^{2}}}\nonumber\\
&=&\lim _{r\rightarrow \infty}|\psi(re^{-i \theta})|^{2}e^{r^{2}(|a|^{2}-1)}\nonumber\nonumber\\
&=&0.
\end{eqnarray}
Since $\psi $ is not identically zero, it is easy to see that there exits a sequence $\{r_{n}\}$ such that for any $n$, $r_{n}$ is a positive real number, $r_{n} \rightarrow \infty$ as $n\rightarrow \infty$ and $\psi(r_{n}e^{-i\theta})\neq 0$
for each $n$. We have $C_{\psi,\varphi}^{\ast}\frac{K_{r_{n}e^{-i \theta}}}{\|K_{r_{n}e^{-i \theta}}\|}=\overline{\psi(r_{n}e^{-i \theta})}\frac{K_{ar_{n}e^{-i \theta}}}{\|K_{r_{n}e^{-i \theta}}\|} \neq 0$ and so $\frac{K_{r_{n}e^{-i \theta}}}{\|K_{r_{n}e^{-i \theta}}\|} \not\in \mbox{Ker}(C_{\psi,\varphi}^{\ast})$. Equation (11) and \cite[Proposition 6.1, p. 363 ]{c1} show that $C_{\psi,\varphi}^{\ast}$ does not have closed range. Thus, $C_{\psi,\varphi}$ does not have closed range (see \cite[Proposition 6.2, p. 364 ]{c1}).
Hence $|a|=1$ and the result follows from the second paragraph of Section 2.\hfill $\Box$ \\ \par

%It shows that $0 \in \sigma_{ap}(C_{\psi,\varphi}^{\ast})$. By  \cite[Proposition 4.4, p. 359]{c1}, $0 \in \sigma_{l,e}(C_{\psi,\varphi}^{\ast}) \cup \{\lambda \in \sigma_{p}(C_{\psi,\varphi}^{\ast}: \mbox{dim ker}(C_{\psi,\varphi}^{\ast}-\lambda)<\infty\}$. If $0 \in  \sigma_{l,e}(C_{\psi,\varphi}^{\ast})$, then by the Open Mapping Theorem and \cite[Proposition 4.3, p. 359]{c1}, $\mbox{ran}(C_{\psi,\varphi}^{\ast})$ is not closed and so $\mbox{ran}(C_{\psi,\varphi})$ is not closed which is a contradiction (see \cite[Proposition 6.2, p. 364]{c1}). Now suppose that $0 \in  \{\lambda \in \sigma_{p}(C_{\psi,\varphi}^{\ast}: \mbox{dim ker}(C_{\psi,\varphi}^{\ast}-\lambda)<\infty\}$ and $0 \notin \sigma_{le}(C_{\psi,\varphi}^{\ast})$. By the Open Mapping Theorem $\mbox{dim ker}(C_{\psi,\varphi})=0$. Since  $0 \notin \sigma_{le}(C_{\psi,\varphi}^{\ast})$, again by \cite[Proposition 4.3, p. 359]{c1}, $C_{\psi,\varphi}^{\ast}$ has closed range. Then \cite[Theorem 2.3 (a), p. 350]{c1} implies that $C_{\psi,\varphi}^{\ast}$ is Fredholm which is a contradiction by Equation (10) and \cite[Theorem 2.3 (c), p. 350]{c1}.
% Hence $|a|=1$ and the result follows from \cite[Proposition 2.1, p. 350]{lef}.\hfill $\Box$ \\ \par

\section{Norm of weighted composition operator}
Suppose that $\varphi$ and $\psi$ are entire functions on $ \mathbb{C}$ and $C_{\psi,\varphi}$ is bounded on $\mathcal{F}^{2}$. It is well-known that for any $w \in \mathbb{C}$,
$$C_{\psi,\varphi}^{\ast}K_{w}=\overline{\psi(w)}K_{\varphi(w)}.$$
We use this formula in the following remark. Furthermore, in this section for each $c \in  \mathbb{C}$, $M_{K_{c}}$ is multiplication by the kernel function $K_{c}$.\\ \par

Remark 4.1.
For $\varphi(z)=az+b$, where $|a| \leq 1$ and $b,c \in \mathbb{C}$, $\|C_{\varphi}\|$ was found in \cite[Theorem 4]{car}.
 In \cite{car}, the inner product of the Fock space is different from ours.
 An analogue of \cite[Theorem 4]{car} holds for $\mathcal{F}^{2}$ with our definition. One can follow the outline of the proof of \cite[Theorem 4]{car} to find $\|C_{\varphi}\|$ on $\mathcal{F}^{2}$. Moreover, we state another proof for finding $\|C_{\varphi}\|$. We break it into two cases. \\
 First assume that $\varphi(z)=az+b$, where $|a|< 1$.
By \cite[Theorem 2]{car}, $C_{\varphi}$ is compact.
Then \cite[Lemma 2]{car} implies that $C_{\varphi}^{\ast}C_{\varphi}=M_{K_{b}}C_{\overline{a}z}C_{az+b}=M_{K_{b}}C_{|a|^{2}z+b}$ is compact ( note that by the similar proof which was stated in \cite[Lemma 2]{car}, we can see that $C_{az+b}^{\ast}=C_{K_{b}, \overline{a}z}$).
 We know that $\|C_{\varphi}\|^{2}=\|C_{\varphi}^{\ast}C_{\varphi}\|=r(M_{K_{b}}C_{|a|^{2}z+b})$. Since $M_{K_{b}}C_{|a|^{2}z+b}$ is compact,
 $r(M_{K_{b}}C_{|a|^{2}z+b})=
\sup \{ |\lambda|: \lambda \in \sigma_{p}(M_{K_{b}}C_{|a|^{2}z+b})\}$ (see \cite[Theorem 7.1,  p. 214 ]{c1}).
 By Proposition 2.1, for each $\lambda \in \sigma_{p}(M_{K_{b}}C_{|a|^{2}z+b})$, $|\lambda|\leq |K_{b}(\frac{b}{1-|a|^{2}})|=e^{\frac{|b|^{2}}{1-|a|^{2}}}$. We have $(M_{K_{b}}C_{|a|^{2}z+b})^{\ast}K_{\frac{b}{1-|a|^{2}}}=e^{\frac{|b|^{2}}{1-|a|^{2}}}K_{\frac{b}{1-|a|^{2}}}$. Then $e^{\frac{|b|^{2}}{1-|a|^{2}}}\in \sigma_{p}((M_{K_{b}}C_{|a|^{2}z+b})^{\ast})$ and so by   \cite[Theorem 7.1,  p. 214 ]{c1},  $e^{\frac{|b|^{2}}{1-|a|^{2}}} \in\sigma_{p} (M_{K_{b}}C_{|a|^{2}z+b})$. Thus, $\|C_{az+b}\|=e^{\frac{1}{2}\frac{|b|^{2}}{1-|a|^{2}}}$. \par
 Now assume that $\varphi(z)=az+b$ with $|a|=1$. By \cite[Theorem 1]{car}, $b=0$. From \cite[Lemma 2]{car}, one can easily see that $\|C_{\varphi}\|^{2}=\|C_{\varphi}^{\ast}C_{\varphi}\|=\|C_{\overline{a}z}C_{az}\|=\|C_{|a|^{2}z}\|=\|C_{z}\|=1$. Then in this case $\|C_{\varphi}\|=1$.\\ \par

 In the preceding sections, we  saw that among weighted composition operators, $C_{K_{c},\varphi}$ is  much important, when   $c \in \mathbb{C}$. Hence in the following theorem, we try to find $\|C_{K_{c},\varphi}\|$.\\ \par

{\bf Theorem 4.2.} {\it Let $\psi=K_{c} $ and $\varphi(z)=az+b$, where $|a| \leq 1$ and $b,c\in \mathbb{C}$. Suppose that $C_{\psi,\varphi}$ is bounded on $\mathcal{F}^{2}$. \\
(a) If $|a|< 1$, then $\|C_{\psi,\varphi}\|=|e^{\overline{c}\frac{b}{1-a}}|e^{\frac{1}{2}\frac{|\frac{c(1-a)}{\overline{a}-1}+b|^{2}}{1-|a|^{2}}} $.\\
(b)  If $|a|= 1$ and $a \neq 1$, then $\|C_{\psi,\varphi}\|=|e^{\frac{|b|^{2}}{1-\overline{a}}}|$.\\
(c) If $a=1$, then $\|C_{\psi,\varphi}\|=e^{\frac{|b|^{2}}{2}}$.}\bigskip

{\bf Proof.} (a) Assume that $|a|< 1$. By the proof of Proposition 2.3, $C_{\psi,\varphi}$ is unitarily equivalent to $\psi(\frac{b}{1-a})C_{az+u(1-a)+b}$,
where $u=\frac{c}{\overline{a}-1}$. Since  $|a|< 1$, \cite[Theorem 4]{car} implies that $\|C_{\psi,\varphi}\|=|\psi(\frac{b}{1-a})|\|C_{az+u(1-a)+b}\|=|\psi(\frac{b}{1-a})| e^{\frac{1}{2}\frac{|\frac{c(1-a)}{\overline{a}-1}+b|^{2}}{1-|a|^{2}}}$ (see also Remark 4.1). \\
(b) Assume that $|a|=1$ and $a \neq 1$. Again by the proof of Proposition 2.3, $\|C_{\psi,\varphi}\|=|\psi(\frac{b}{1-a})|\|C_{az+\frac{c(1-a)}{\overline{a}-1}+b}\|=|\psi(\frac{b}{1-a})|\|C_{az+ca+b}\|$. Since $|a|=1$ and $C_{K_{c},\varphi}$ is bounded, from  the second paragraph of  Section 2, $c=-\overline{a}b$.  Then $ca+b=0$. Therefore, $\|C_{\psi,\varphi}\|=|\psi(\frac{b}{1-a})|\|C_{az}\|=|\psi(\frac{b}{1-a})|$ (see \cite[Theorem 4]{car} and  Remark 4.1). Thus, $\|C_{\psi,\varphi}\|=|e^{\frac{-a|b|^{2}}{1-a}}|=|e^{\frac{|b|^{2}}{1-\overline{a}}}|$.\\
 (c) Assume that $a=1$. As we saw in the second paragraph of Section 2, $c=-\overline{a}b=-b$. Then $\psi(z)=e^{-\overline{b}z}$.
 Now we must find $\|C_{e^{-\overline{b}z}, z+b}\|$. By \cite[Corollary 1.2]{zhao1}, $e^{\frac{-|b|^{2}}{2}}C_{e^{-\overline{b}z},z+b}$ is unitary. Hence $\|C_{\psi,\varphi}\|=e^{\frac{|b|^{2}}{2}}$.\hfill $\Box$ \\ \par

%\textbf{Acknowledgments}
%\bigskip

\bigskip
{M. Fatehi, Department of Mathematics, Shiraz Branch, Islamic Azad
University, Shiraz, Iran. \par E-mail: fatehimahsa@yahoo.com \par
}

\end{document}